\pdfoutput=1

\title{Genus expansion for real Wishart matrices}
\author{C. Emily I. Redelmeier}

\documentclass[11pt]{article}

\usepackage{graphicx}
\usepackage{graphics}
\usepackage{amsfonts}
\usepackage{amssymb}
\usepackage{amsthm}
\usepackage{amsmath}
\usepackage{color}
\usepackage{cite}

\newtheorem{theorem}{Theorem}[section]
\newtheorem{lemma}[theorem]{Lemma}
\newtheorem{proposition}[theorem]{Proposition}

\theoremstyle{remark}
\newtheorem{remark}[theorem]{Remark}
\newtheorem{notation}[theorem]{Notation}
\newtheorem{example}[theorem]{Example}

\theoremstyle{definition}
\newtheorem{definition}[theorem]{Definition}

\begin{document}

\maketitle

\begin{abstract}
We present an exact formula for moments and cumulants of several real compound Wishart matrices in terms of an Euler characteristic expansion, similar to the genus expansion for complex random matrices.  We consider their asymptotic values in the large matrix limit: as in a genus expansion, the terms which survive in the large matrix limit are those with the greatest Euler characteristic, that is, either spheres or collections of spheres.  This topological construction motivates an algebraic expression for the moments and cumulants in terms of the symmetric group.  We examine the combinatorial properties distinguishing the leading order terms.  By considering higher cumulants, we give a central limit-type theorem for the asymptotic distribution around the expected value.
\end{abstract}

\section{Introduction}

Wishart matrices are random matrices of the form \(AX^{T}BXC\), where \(X\) is a matrix with Gaussian entries and \(A\), \(B\) and \(C\) are constant matrices which may be required to be symmetric, positive definite, or the identity matrix.  This matrix ensemble was introduced in Wishart's paper of 1928 \cite{Wishart} in order to study the expected covariances of multivariate data.  The moments of Wishart matrices, as well as more general symmetric polynomials in their eigenvalues, are studied in such papers as \cite{HSS, MR2132270, MR2439565, MR2240781, MR2337139}.  The behaviour of these quantities in large matrices and the connections to free probability are considered in papers such as \cite{MR2052516, MR2240781, MR2337139}.

In this paper, we exhibit an exact formula for the expected value of products of traces of powers of Wishart matrices in (\ref{main formula}):
\begin{multline*}
\mathbb{E}\left(\mathrm{Tr}\left(X^{T}D_{1}XD_{2}\cdots X^{T}D_{2n_{1}-1}XD_{2n_{1}}\right)\cdots\right.\\\left.\mathrm{Tr}\left(X^{T}D_{2n-2n_{r}+1}XD_{2n-2n_{r}+2}\cdots X^{T}D_{2n-1}XD_{2n}\right)\right)\\=\sum_{\pi}\mathrm{Tr}_{\gamma_{-}^{-1}\hat{\pi}\gamma_{+}/2}\left(D_{1},\ldots,D_{2n}\right)
\end{multline*}
where the sum on the right hand side is over all pairings \(\pi\) on \(2n\) elements, interpreted as permutations, \(\hat{\pi}\) is a permutation derived from \(\pi\), and \(\gamma_{-}\) and \(\gamma_{+}\) are fixed permutations dependent on the left hand side expression.  If we wish to calculate cumulants, rather than moments, we sum only over \(\pi\) such that the subgroup \(\langle \gamma_{+}\gamma_{-}^{-1},\hat{\pi}\rangle\) is transitive.

If the permutations are interpreted as the information to construct a topological cell-complex, this expression becomes a genus expansion similar to those for complex matrix ensembles (see \cite{MR1492512, MR2036721, MR2052516}), but which also includes terms corresponding to non-orientable surfaces.  As in the complex cases, the order of the terms depends only on the Euler characteristic of the corresponding surface, allowing us to distinguish which terms survive in the asymptotic limit.  The leading terms are again those with the highest Euler characteristic.  The terms contributing to the mixed moments of traces are arbitrary surfaces, so the highest order terms are those whose corresponding surface consists of as many disjoint spheres as possible.  The terms contributing to cumulants of the traces are connected surfaces, so the leading terms are those corresponding to spheres.

This topological interpretation is also conducive to studying several independent Wishart matrices, or the more general case of Wishart matrices indexed by a Hilbert space whose covariances are given by the inner product.  This construction also allows us to handle expressions in which the random matrix \(X\) and its transpose appear arbitrarily in the products of matrices, rather than alternatingly, as in the usual Wishart matrix model.

The topological interpretation motivates the algebraic expression for these expected values, following the expression in \cite{MR2052516} for the complex case based on the interpretation of a map or hypermap as a pair of permutations in \cite{MR0404045}.  We now must consider non-orientable surfaces as well as orientable ones.  The algebraic formulas in \cite{MR0404045} require a consistent orientation on the surface, so we construct orientable two-sheeted covering spaces of our original surfaces and extend the algebraic formula to this situation.  

In Section~\ref{notation}, we will outline the notation and results we will be using and define the matrix model which we will be considering.  To simplify our notation, we will delay discussion of the several matrix model and expressions with arbitrary transposes until Section~\ref{generalisations}; however, most of the arguments are identical.  In Section~\ref{geometry}, we provide the geometric motivation for the algebraic formulas using a construction similar to \cite{MR1492512, MR2036721}.  For a rigorous proof, we prefer to use an algebraic construction for the moments and cumulants using the representations for graphs in surfaces developed in \cite{MR0404045} and used in the complex Wishart case in \cite{MR2052516}, which we describe in Section~\ref{algebra}.  The proof of the main formula is given in Section~\ref{proof}, and its asymptotic behaviour is discussed in Section~\ref{asymptotics}.

\section{Notation, Definitions and Lemmas}
\label{notation}

\begin{notation}
Denote the set \(\left\{1,\ldots,n\right\}\) by \(\left[n\right]\) and the set \(\left\{1,\ldots,n\right\}\cup\left\{-1,\ldots,-n\right\}\) by \(\left[\pm n\right]\).
\end{notation}

We write the usual trace \(\mathrm{Tr}\) and the normalised trace \(\mathrm{tr}:=\frac{1}{N}\mathrm{Tr}\).

\subsection{The Matrix Model}

\begin{definition}
\label{Wishart matrix}
We will denote by \(X\) an \(M\times N\) matrix whose \(ij\)th entry is \(\frac{1}{\sqrt{N}}f_{ij}\), where the \(f_{ij}\) are independent real \(N\left(0,1\right)\) random variables (that is, centred Gaussians with standard deviation \(1\)).

A matrix of the form \(X^{T}D_{k}XD_{k+1}\), where \(D_{k}\) is an arbitrary constant \(M\times M\) matrix and \(D_{k+1}\) is an arbitrary \(N\times N\) matrix, is called a {\em real Wishart matrix}.
\end{definition}

We will be considering the statistics of these matrices as \(N\rightarrow\infty\), and we assume \(M\rightarrow\infty\) at the same time.  We will continue to consider the normalised trace \(\frac{1}{N}\mathrm{Tr}\), even when we are considering an \(M\times M\) matrix.  We will also assume that the large \(N\) limits of the normalised traces of products of \(D_{k}\) matrices approach a finite limit.

Let \(n_{1},\ldots,n_{r}\) be positive integers.  We will be considering the moments and cumulants of the traces of products of Wishart matrices, that is random variables of the form
\begin{multline}
Y_{k}:=\mathrm{tr}\left(X^{T}D_{2\left(n_{1}+\cdots+n_{k-1}\right)+1}XD_{2\left(n_{1}+\cdots+n_{k-1}\right)+2}\cdots\right.\\\left.X^{T}D_{2\left(n_{1}+\cdots+n_{k}\right)-1}XD_{2\left(n_{1}+\cdots+n_{k}\right)}\right)\textrm{.}\label{trace random variable}
\end{multline}

There are no requirements on the \(D_{k}\) matrices: they need not be symmetric, and none are required to be equal.  Both the subscripts of the matrices \(D_{k}\) and the appearances of transposes appear naturally in our formula.  We will denote the transpose of \(D_{k}\) by \(D_{-k}\) and the \(ij\)th entry of \(D_{k}\) by \(d^{\left(k\right)}_{ij}\).

Since we only consider the traces of products of matrices, expressions which begin with a constant matrix can be calculated by cycling any initial matrix to the end of the expression inside the trace.

\subsection{Partitions}

\begin{definition}
A {\em partition} of a set \(S\) is a set of subsets \(V_{1},\ldots,V_{m}\subseteq S\) called {\em blocks} such that \(V_{i}\neq\emptyset\) for \(1\leq i\leq m\), \(V_{i}\cap V_{j}=\emptyset\) for all \(i\neq j\), and \(V_{1}\cup\cdots\cup V_{m}=S\).  We denote the set of all partitions of \(S\) by \({\cal P}\left(S\right)\), and the set of all partitions of \(\left[n\right]\) by \({\cal P}\left(n\right)\).

A partitions whose blocks contain exactly two elements is called a {\em pairing}.  We denote the set of all pairings of a set \(S\) by \({\cal P}_{2}\left(S\right)\), and the set of all pairings of \(\left[n\right]\) by \({\cal P}_{2}\left(n\right)\).
\end{definition}

\subsection{Moments and Cumulants}

\begin{definition}
The \(n\)th mixed moment of random variables \(X_{1},\ldots,X_{n}\) is defined as
\[a_{n}\left(X_{1},\ldots,X_{n}\right):=\mathbb{E}\left(X_{1}\cdots X_{n}\right)\textrm{.}\]
\end{definition}

\begin{definition}
Cumulants are functions \(k_{1},k_{2},\ldots\) defined so that \(k_{n}\) is an \(n\)-linear function, and satisfying the moment-cumulant formula:
\begin{equation}
\label{moment-cumulant formula}a_{n}\left(X_{1},\ldots,X_{n}\right)=\sum_{\pi\in{\cal P}\left(n\right)}\prod_{V=\left\{i_{1},\ldots,i_{m}\right\}\in\pi}k_{m}\left(X_{i_{1}},\ldots,X_{i_{m}}\right)
\end{equation}
for all random variables \(X_{1},\ldots,X_{n}\) and all \(n\geq 0\).
\end{definition}

It can easily be shown by induction that for all \(n\), the moment cumulant formulas up to \(n\) along with the cumulants \(k_{1},\ldots,k_{n-1}\) uniquely define the multilinear function \(k_{n}\).  The first cumulant is the expected value, and the second is the covariance.

We denote the term on the left hand side of the moment-cumulant formula corresponding to a given partition \(\pi\in{\cal P}\left(n\right)\) by
\[k_{\pi}\left(X_{1},\ldots,X_{n}\right):=\prod_{V=\left\{i_{1},\ldots,i_{m}\right\}\in\pi}k_{m}\left(X_{i_{1}},\ldots,X_{i_{m}}\right)\textrm{.}\]

\subsection{The Wick Formula}

The following formula allows us to express the expectations of products of Gaussian random variables in terms of a combinatorial construction.

\begin{lemma}[The Wick Formula]
\label{Wick formula}
Let \(\left\{f_{\lambda}:\lambda\in\Lambda\right\}\) be a centred family of Gaussian random variables (that is, each \(f_{\lambda}\) is a linear combination of random variables from a collection of independent \(N\left(0,1\right)\) random variables).  Then
\[\mathbb{E}\left(f_{\lambda_{1}}\cdots f_{\lambda_{n}}\right)=\sum_{\pi\in{\cal P}_{2}\left(n\right)}\prod_{\left\{k,l\right\}\in\pi}\mathbb{E}\left(f_{\lambda_{k}}f_{\lambda_{l}}\right)\textrm{.}\]
\end{lemma}
In particular, if \(n\) is odd, there are no pairings, and the expected value is \(0\).

For a proof, see \cite{MR1474726}, theorem 1.23.

\subsection{Traces Along Permutations}

For a permutation \(\pi\in S_{n}\) with cycle notation \((c_{1},\ldots,c_{n_{1}})\cdots(c_{n_{1}+\ldots+n_{r-1}+1},\allowbreak\ldots,c_{n})\), we define the trace along \(\pi\) of matrices \(A_{1},\ldots,A_{n}\) as
\[\mathrm{Tr}_{\pi}\left(A_{1},\ldots,A_{n}\right):=\mathrm{Tr}(A_{c_{1}}\cdots A_{c_{n_{1}}})\cdots\mathrm{Tr}(A_{c_{n_{1}+\cdots+n_{r-1}+1}}\cdots A_{c_{n}})\textrm{.}\]
Because the trace is cyclic, this expression is well defined.  For permutations on a set of signed integers, we let a negative integer denote a transpose: if we let \(A_{k}^{\left(1\right)}=A_{k}\) while \(A_{k}^{\left(-1\right)}=A_{k}^{T}\), and if \(\pi=(c_{1},\ldots,c_{n_{1}})\cdots(c_{n_{1}+\ldots+n_{r-1}+1}\cdots\allowbreak c_{n})\) is a permutation of a subset of \(\left[\pm m\right]\), then
\begin{multline*}
\mathrm{Tr}_{\pi}\left(A_{1},\ldots,A_{m}\right):=\mathrm{Tr}\left(A_{\left|c_{1}\right|}^{\left(\mathrm{sgn}\left(c_{1}\right)\right)}\cdots A_{\left|c_{n_{1}}\right|}^{\left(\mathrm{sgn}\left(c_{n_{1}}\right)\right)}\right)\cdots\\\mathrm{Tr}\left(A_{\left|c_{n_{1}+\cdots+n_{r-1}+1}\right|}^{\left(\mathrm{sgn}\left(c_{n_{1}+\cdots+n_{r-1}+1}\right)\right)}\cdots A_{\left|c_{n}\right|}^{\left(\mathrm{sgn}\left(c_{n}\right)\right)}\right)\textrm{.}
\end{multline*}

We will be making use of the following folklore result:
\begin{lemma}
\label{trace along permutation}
Let \(A_{1},\ldots,A_{n}\) be matrices such that \(A_{k}\) is an \(N_{k}\times N_{\sigma\left(k\right)}\) matrix.  If we denote the \(ij\)th entry of matrix \(A_{k}\) by \(a^{\left(k\right)}_{ij}\), then
\[\mathrm{Tr}_{\sigma}\left(A_{1},\ldots,A_{n}\right)=\sum_{\substack{1\leq i_{1}\leq N_{1}\\\vdots\\1\leq i_{n}\leq N_{n}}}a^{\left(1\right)}_{i_{1}i_{\sigma\left(1\right)}}\cdots a^{\left(n\right)}_{i_{n}i_{\sigma\left(n\right)}}\textrm{.}\]
\end{lemma}
This formula may be proven by standard calculation.  The subscripts on the matrices \(A_{k}\) may be replaced by the elements of a finite index set and \(\sigma\) by a permutation of that set, allowing us to use this formula in the situation where the index set is a subset of \(\left[\pm n\right]\).

\section{Geometric Motivation}
\label{geometry}

In this section we sketch a method of calculation using surface gluings.  We delay a rigorous proof until Section~\ref{proof}, which will use only Lemma~\ref{lemma} from this section.

We begin by examining the moments of the random variables \(Y_{k}\) defined in (\ref{trace random variable}).  We will let \(n=n_{1}+\cdots+n_{r}\) (\(n_{k}\) as in (\ref{trace random variable})).

\begin{notation}
\label{gamma}
Let \(\gamma:=\left(1,\ldots,2n_{1}\right)\cdots\left(2\left(n_{1}+\cdots+n_{r-1}\right)+1,\ldots,2n\right)\in S_{2n}\).
\end{notation}

\begin{lemma}
\label{lemma}
A mixed moment of the \(Y_{k}\) may be expressed in terms of the entries of the matrices:
\begin{multline}
a_{n}\left(Y_{1},\ldots,Y_{r}\right)\\=N^{-n-r}\sum_{\pi\in{\cal P}_{2}}\sum_{\substack{1\leq i_{1},\ldots,i_{2n}\leq N\\1\leq j_{1},\ldots,j_{2n}\leq M}}d^{\left(1\right)}_{i_{1}i_{\gamma\left(1\right)}}d^{\left(2\right)}_{j_{2}j_{\gamma\left(2\right)}}\cdots d^{\left(2n-1\right)}_{i_{2n-1}i_{\gamma\left(2n-1\right)}}d^{\left(2n\right)}_{j_{2n}j_{\gamma\left(2n\right)}}\\\prod_{\left\{k,l\right\}\in\pi}\mathbb{E}\left(f_{i_{k}j_{k}}f_{i_{l}j_{l}}\right)\label{sum over pairings}\textrm{.}
\end{multline}
Furthermore, since we can express the product of expected values as:
\begin{equation}
\prod_{\left\{k,l\right\}\in\pi}\mathbb{E}\left(f_{i_{k}j_{k}}f_{i_{l}j_{l}}\right)=\left\{\begin{array}{ll}1,&\textrm{\(i_{k}=i_{l}\) and \(j_{k}=j_{l}\) for all \(\left\{k,l\right\}\in\pi\)}\\0,&\textrm{otherwise}\end{array}\right.\label{expected value expression}
\end{equation}
we may replace the expected value expression with a condition on the indices:
\begin{multline*}
a_{r}\left(Y_{1},\ldots,Y_{r}\right)=N^{-n-r}\sum_{\pi\in{\cal P}_{2}\left(2n\right)}\sum_{\substack{i:\left[2n\right]\rightarrow\left[N\right]:i=i\circ\pi\nonumber\\j:\left[2n\right]\rightarrow\left[M\right]:j=j\circ\pi}}d^{\left(1\right)}_{i_{1}i_{\gamma\left(1\right)}}d^{\left(2\right)}_{j_{2}j_{\gamma\left(2\right)}}\cdots\\d^{\left(2n-1\right)}_{i_{2n-1}i_{\gamma\left(2n-1\right)}}d^{\left(2n\right)}_{j_{2n}j_{\gamma\left(2n\right)}}
\end{multline*}
where \(i\left(k\right)=i_{k}\) and \(j\left(k\right)=j_{k}\).

\begin{proof}
We can write the moment:
\begin{eqnarray*}
a_{r}\left(Y_{1},\cdots,Y_{r}\right)&=&\mathbb{E}\left(Y_{1}\cdots Y_{r}\right)
\\&=&\mathbb{E}\left(\mathrm{tr}\left(X^{T}D_{1}XD_{2}\cdots X^{T}D_{2n_{1}-1}XD_{2n_{1}}\right)\cdots\right.\\&&\left.\mathrm{tr}\left(X^{T}D_{2\left(n_{1}+\cdots+n_{r-1}\right)+1}XD_{2\left(n_{1}+\cdots n_{r-1}\right)+2}\cdots\right.\right.\\&&\left.\left.X^{T}D_{2n-1}XD_{2n}\right)\right)\textrm{.}
\end{eqnarray*}
We can interpret this expression as a trace along permutation \(\gamma\), so according to Lemma~\ref{trace along permutation},
\begin{eqnarray*}
a_{r}\left(Y_{1},\cdots,Y_{r}\right)&=&\mathbb{E}\left(N^{-r}\mathrm{Tr}_{\gamma}\left(X^{T}D_{1},XD_{2},\ldots,X^{T}D_{2n-1},X,D_{2n}\right)\right)
\\&=&N^{-r}\sum_{\substack{1\leq i_{1},\ldots,i_{2n}\leq N\\1\leq j_{1},\ldots,j_{2n}\leq M}}\mathbb{E}\left(X^{T}_{j_{1}i_{1}}d^{\left(1\right)}_{i_{1}i_{\gamma\left(1\right)}}X_{i_{2}j_{2}}d^{\left(2\right)}_{j_{2}j_{\gamma\left(2\right)}}\cdots\right.\\&&\left.X^{T}_{i_{2n-1}j_{2n-1}}d^{\left(2n-1\right)}_{i_{2n-1}i_{\gamma\left(2n-1\right)}}X_{i_{2n}j_{2n}}d^{\left(2n\right)}_{j_{2n}j_{\gamma\left(2n\right)}}\right)
\\&=&N^{-n-r}\sum_{\substack{1\leq i_{1},\ldots,i_{2n}\leq N\\1\leq j_{1},\ldots,j_{2n}\leq M}}d^{\left(1\right)}_{i_{1}i_{\gamma\left(1\right)}}d^{\left(2\right)}_{j_{2}j_{\gamma\left(2\right)}}\cdots\\&&d^{\left(2n-1\right)}_{i_{2n-1}i_{\gamma\left(2n-1\right)}}d^{\left(2n\right)}_{j_{2n}j_{\gamma\left(2n\right)}}\mathbb{E}\left(f_{i_{1}j_{1}}\cdots f_{i_{2n}j_{2n}}\right)\textrm{.}
\end{eqnarray*}
Applying the Wick formula (Lemma~\ref{Wick formula}) to the expected value expression in the above, we get:
\begin{eqnarray*}
a_{r}\left(Y_{1},\ldots,Y_{r}\right)&=&N^{-n-r}\sum_{\substack{1\leq i_{1},\ldots,i_{2n}\leq N\\1\leq j_{1},\ldots,j_{2n}\leq M}}d^{\left(1\right)}_{i_{1}i_{\gamma\left(1\right)}}d^{\left(2\right)}_{j_{2}j_{\gamma\left(2\right)}}\cdots\\&&d^{\left(2n-1\right)}_{i_{2n-1}i_{\gamma\left(2n-1\right)}}d^{\left(2n\right)}_{j_{2n}j_{\gamma\left(2n\right)}}\sum_{\pi\in{\cal P}_{2}}\prod_{\left\{k,l\right\}\in\pi}\mathbb{E}\left(f_{i_{k}j_{k}}f_{i_{l}j_{l}}\right)
\\&=&N^{-n-r}\sum_{\pi\in{\cal P}_{2}}\sum_{\substack{1\leq i_{1},\ldots,i_{2n}\leq N\\1\leq j_{1},\ldots,j_{2n}\leq M}}d^{\left(1\right)}_{i_{1}i_{\gamma\left(1\right)}}d^{\left(2\right)}_{j_{2}j_{\gamma\left(2\right)}}\cdots\\&&d^{\left(2n-1\right)}_{i_{2n-1}i_{\gamma\left(2n-1\right)}}d^{\left(2n\right)}_{j_{2n}j_{\gamma\left(2n\right)}}\prod_{\left\{k,l\right\}\in\pi}\mathbb{E}\left(f_{i_{k}j_{k}}f_{i_{l}j_{l}}\right)
\end{eqnarray*}
where we have reversed the order of summation to express the moment as a sum over pairings.  This proves (\ref{sum over pairings}).

In the term corresponding to a given pairing \(\pi\), the factor \(\mathbb{E}\left(f_{i_{k}j_{k}}f_{i_{l}j_{l}}\right)\), and therefore the entire term, will be equal to zero if \(f_{i_{k}j_{k}}\) is independent from \(f_{i_{l}j_{l}}\).  So we find:
\[\prod_{\left\{k,l\right\}\in\pi}\mathbb{E}\left(f_{i_{k}j_{k}}f_{i_{l}j_{l}}\right)=\left\{\begin{array}{ll}1,&\textrm{\(i_{k}=i_{l}\) and \(j_{k}=j_{l}\) for all \(\left\{k,l\right\}\in\pi\)}\\0,&\textrm{otherwise}\end{array}\right.\]
proving (\ref{expected value expression}).  We can interpret the conditions for a nonzero expected value as a further set of constraints on the indices:
\begin{multline*}
a_{r}\left(Y_{1},\ldots,Y_{r}\right)=N^{-n-r}\sum_{\pi\in{\cal P}_{2}\left(2n\right)}\sum_{\substack{i:\left[2n\right]\rightarrow\left[N\right]:i=i\circ\pi\nonumber\\j:\left[2n\right]\rightarrow\left[M\right]:j=j\circ\pi}}d^{\left(1\right)}_{i_{1}i_{\gamma\left(1\right)}}d^{\left(2\right)}_{j_{2}j_{\gamma\left(2\right)}}\cdots\\d^{\left(2n-1\right)}_{i_{2n-1}i_{\gamma\left(2n-1\right)}}d^{\left(2n\right)}_{j_{2n}j_{\gamma\left(2n\right)}}\textrm{.}
\end{multline*}
This completes the proof of the lemma.
\end{proof}
\end{lemma}

\begin{figure}
\centering
\scalebox{0.75}{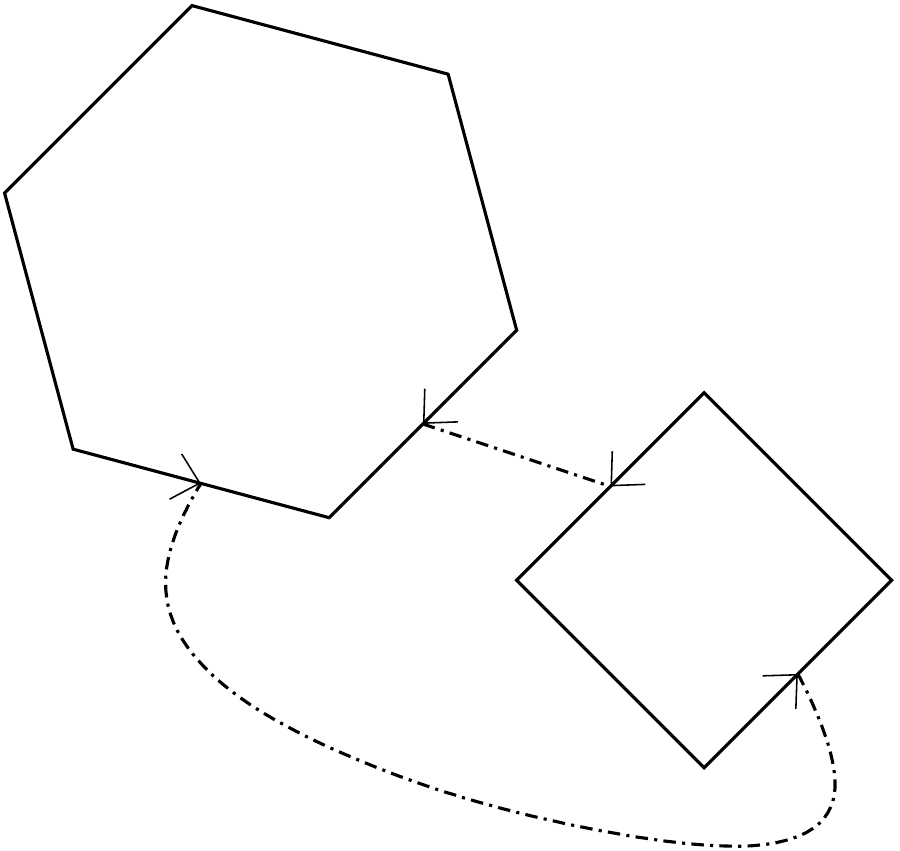}
\caption{The faces constructed in the calculation of \(\mathbb{E}\left(\mathrm{tr}\left(X^{T}D_{1}XD_{2}\cdots X^{T}D_{5}XD_{6}\right)\mathrm{tr}\left(X^{T}D_{7}XD_{8}X^{T}D_{9}XD_{10}\right)\right)\), with the random variables \(f_{i_{2}j_{2}}\) and \(f_{i_{7}j_{7}}\) (from \(X\) and \(X^{T}\) respectively) paired with each other, and \(f_{i_{1}j_{1}}\) and \(f_{i_{9}j_{9}}\) (both from \(X^{T}\)) paired with each other.  The resulting edge identifications are shown.}
\label{faces}
\end{figure}

These constraints can be expressed geometrically as follows.  For each random variable \(Y_{k}\), we construct a \(2n_{k}\)-gon.  We associate each edge with a random matrix term and each vertex with a constant matrix term, cyclically (counterclockwise) in the order they appear in the expression for \(Y_{k}\).  We label each end of each edge with the index the edge and the vertex at that end share.  An example is shown in Figure~\ref{faces}.

A pairing \(\pi\) constrains \(i_{k}=i_{l}\) and \(j_{k}=j_{l}\) for each \(\left\{k,l\right\}\in\pi\).  We represent this by identifying the edges \(X_{k}\) and \(X_{l}\) so that the edge-ends which have been labelled \(i_{k}\) and \(i_{l}\) line up, as do those labelled \(j_{k}\) and \(j_{l}\).  Figure~\ref{faces} shows two examples.  A pairing of a random variable from an \(X\) term with one from an \(X^{T}\) term will result in an untwisted edge identification while a pairing of random variables from two \(X\) terms or two \(X^{T}\) terms will result in a twisted pairing.  The resulting surface thus may not be orientable.  However, since each edge is paired with exactly one other, it will be a compact surface without boundary.

\begin{figure}
\centering
\scalebox{0.75}{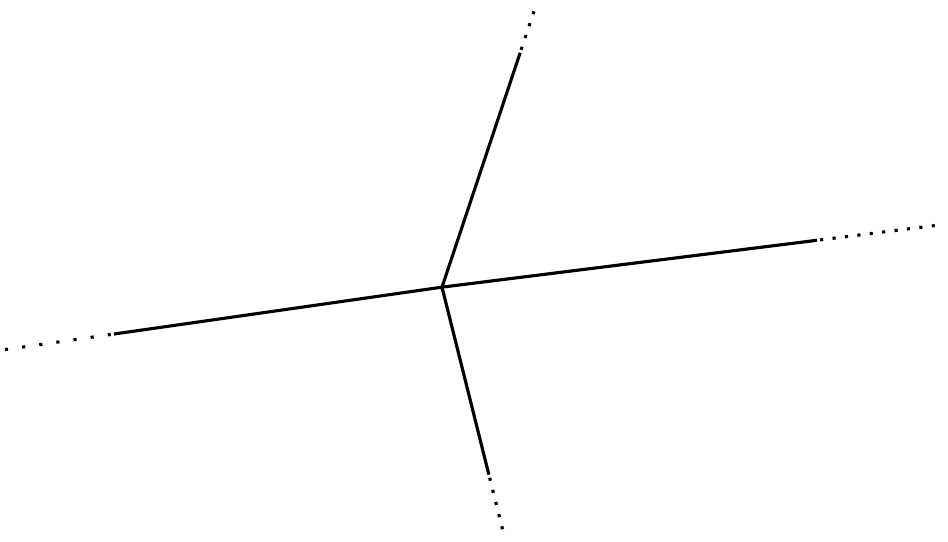}
\caption{A vertex that appears in the surface resulting from the pairing \(\left\{\left\{1,9\right\}\left\{2,7\right\}\left\{3,4\right\}\left\{5,10\right\}\left\{6,8\right\}\right\}\).  To perform the twisted edge identifications, the corners of faces containing the matrices \(D_{5}\) and \(D_{9}\) were flipped over, so the matrices \(D_{5}^{T}\) and \(D_{9}^{T}\) appear instead, since it is \(D_{5}^{T}\) whose \(i_{6}i_{5}\)th entry appears, and likewise with \(D_{9}^{T}\).  The indices are identified as they would be in the expression for the trace \(\mathrm{Tr}\left(D_{1}D_{7}D_{5}^{T}D_{9}^{T}\right)\).}
\label{vertex}
\end{figure}

To determine the constraints on the indices of the constant matrices, we examine the vertices of the resulting surface.  Figure~\ref{vertex} shows a an example.  Due to the presence of twisted edge identifications, some corners may appear flipped over relative to others.  We replace their matrices with their transposes so the indices of the matrix as written always appear in counterclockwise order.

We note that around a vertex the indices of the matrices are identified as they would be in the expression for a trace.  When we express the moment in terms of normalised traces, each trace gives us a factor of \(N\).  The order of the term corresponding to a given pairing whose surface has \(\#V\) vertices will be \(N^{\#V-n-r}\), or, since the number of edges (after identifications) is \(n\) and the number of faces is \(r\), \(N^{\chi-2r}\), where \(\chi=\#V-\#E+\#F\) is the Euler characteristic of the surface.  The highest order terms will be those with the greatest Euler characteristic.

\begin{remark}
We note that this construction appears implicitly in many of the apparently different approaches to the calculation of moments.

The face gluings we have constructed can be interpreted as the duals of the M\"{o}bius graphs in \cite{MR2480549, MR2187951}.

In the constructions used in \cite{HSS, MR2132270, MR2439565}, two graphs are constructed whose union consists of cycles whose lengths are \(2n_{1},\ldots,2n_{r}\).  We can interpret these cycles as the faces in our construction, although the vertices are edges and the edges are vertices.  The moment is then expressed as a sum over all pairings of the vertices, equivalent to our sum over edge identifications.  The cycles of the unions of the pairing with each of the first two graphs can be interpreted as the vertices in our construction.  The two graphs contributing to the faces keeps track of the orientation on what may be nonorientable surfaces.  However, these constructions cannot easily be extended to the case where the transposes appear arbitrarily.

The closed words in \cite{MR2222385} may be interpreted as faces, and the condition that the graph of a word be a tree when there is one word and a graph with exactly one (possibly degenerate) cycle when there are more words can be translated into a condition on the number of vertices and hence on the Euler characteristic of the resulting surface.

Approaches such as \cite{MR2052516, MR2240781, MR2337139} may be interpreted in terms of the algebraic construction in Section~\ref{algebra}.
\end{remark}

Any compact, connected surface without boundary is either a sphere, the connected sum of some number of tori, or the connected sum of some number of projective planes (for a proof of this classification theorem see, for example, \cite{MR1280460}).  The sphere has Euler characteristic \(\chi=2\), the connected sum of \(n\) tori has \(\chi=2-2n\), and the connected sum of \(n\) projective planes has \(\chi=2-n\).  Our surface will have the highest Euler characteristic, and hence the highest order in \(N\), if it consists of as many spheres as possible.  The leading terms therefore consist of a disjoint sphere for each face, with Euler characteristic \(2r\), giving it order \(N^{0}=1\).  These terms correspond to disjoint noncrossing pairings on each face.

We will see (Proposition~\ref{cumulants}) that the cumulants are the sums of the terms corresponding to connected surfaces.  (The connected components of any surface partition the set of faces, inducing a partition of \(\left[r\right]\).  For any such partition, the contribution to the moment is the product of the contribution of each of its components, so the contributions of connected components satisfy the moment-cumulant formula (\ref{moment-cumulant formula}).)  The highest order terms contributing to the cumulant \(k_{r}\left(Y_{1},\ldots,Y_{r}\right)\) correspond to spheres, of order \(N^{2-2r}\).

In the large \(N\) limit, only the first cumulant survives.  We can also consider the asymptotic distribution of the \(Y_{k}\) around their means by examining their higher cumulants.  The second cumulant is of higher order in \(N\) than the higher cumulants, so asymptotically the distribution is Gaussian (see, for example, \cite{MR1474726}).  The second cumulants can be interpreted as the asymptotic covariances of the \(Y_{k}-\mathbb{E}\left(Y_{k}\right)\).

\section{An Algebraic Expression for the Moments}
\label{algebra}

Following \cite{MR2052516}, we now use the methods of \cite{MR0404045} to algebraically calculate face information from vertex information and edge information (or, in our case the dual problem of calculating vertex information from face information and edge information).  Each end of each edge is labelled.  We then construct a permutation \(\gamma\) of these labels which permutes the labels of the edge-ends incident to each vertex in a counter-clockwise direction, and a permutation \(\pi\) which swaps the labels of the edge-ends of each edge.  The cycles of the permutation \(\pi\circ\gamma\) then enumerate an edge-end for each edge appearing around each face.  In this case, however, our surfaces may be nonorientable, so we will consider an orientable two-sheeted covering space of each of our surfaces.

An orientable two-sheeted covering space may be defined for any surface (see, for example, \cite{MR1867354}, pages 234--235).  This surface is the one experienced by someone on the original surface, rather than within it.  This construction can be performed more simply in our situation, where we have constructed our surface as a cell complex.  For each face, we construct two preimage faces, corresponding to the two orientations, which we will think of as the front and the back.  For each untwisted edge identification in our original surface, we identify the preimage edges of the front faces in the same manner, and likewise the preimage edges of the back faces.  For each twisted edge identification, we identify the preimage edge of each front face with the preimage edge of each back face.

To accommodate this two-sheeted covering, we double the set \(\left[2n\right]\) on which our permutations act.  We will use negative integers to represent the ``back'' preimage points, so our set will now be \(\left[\pm 2n\right]\).  We construct a pairing \(\hat{\pi}\) which connects the front and back appropriately.  We will do this by connecting \(k\) with \(l\) and \(-k\) with \(-l\) for every pair \(\left\{k,l\right\}\in\pi\), then ``twisting'' each pairing by pairing \(k\) with \(-l\) and \(-k\) with \(l\).  We will then twist each end which corresponds to an \(X^{T}\) term.  This will result in a total of \(1\) twist if both terms are \(X\) terms, \(2\) twists (equivalent to no twists) if an \(X\) term is paired with an \(X^{T}\) term, and \(3\) twists (equivalent to \(1\) twist) if two \(X^{T}\) terms are paired.

\begin{notation}
\label{delta}
Let \(\delta\) be the permutation on \(\left[\pm 2n\right]\) which reverses the sign on its argument.

Let \(\delta^{\prime}\) be the permutation that takes \(k\) to \(-k\) and \(-k\) to \(k\) if the \(k\)th appearance of \(X\) in our original expression is transposed and does nothing otherwise.
\end{notation}

Conjugating \(\pi\) by \(\delta\) constructs a copy of the pairing on the negative integers, and the pairing \(\delta\pi\delta\pi\) consists of both copies of the pairing.  To put a twist in each pairing, we multiply by \(\delta\) again, giving us \(\delta\delta\pi\delta\pi=\pi\delta\pi\).  We now wish to twist each end corresponding to an \(X^{T}\).  We do so by conjugating by \(\delta^{\prime}\)

\begin{notation}
\label{pi hat}
Let \(\hat{\pi}:=\delta^{\prime}\pi\delta\pi\delta^{\prime}\).
\end{notation}

We construct permutations taking the place of \(\gamma\), one acting on the fronts (positive integers) and the other acting on the backs (negative integers).

\begin{notation}
\label{more gammas}
Let \(\gamma_{+}:=\left(1,\ldots,2n_{1}\right)\cdots\left(2\left(n_{1}+\cdots+n_{r-1}\right)+1,\ldots,2n\right)\) and \(\gamma_{-}:=\left(-1,\ldots,-2n_{1}\right)\cdots\left(2\left(n_{1}+\cdots+n_{r-1}\right)-1,\ldots,-2n\right)\).
\end{notation}

We would expect the Kreweras complement to be \(\hat{\pi}\gamma_{+}\gamma_{-}^{-1}\); however, because of the conventions we are using for numbering the edges and corners, we find that instead of negative integer \(-k\), we get \(\gamma_{-}\left(-k\right)\) (while we find the correct indices when the integers are positive).  To correct this, we conjugate by \(\gamma_{-}^{-1}\), leaving us with the permutation \(\gamma_{-}^{-1}\hat{\pi}\gamma_{+}\).

The cycles of this permutation enumerate the subscripts of the matrices appearing around each vertex, negative if the matrix must be transposed.  However, since we are using a two-sheeted covering space of the surface, each vertex appears twice.  Our application will require us to discard the redundant cycles.

\begin{example}
Let \(\pi=\left(1,9\right)\left(2,7\right)\left(3,4\right)\left(5,10\right)\left(6,8\right)\).  Then if \(\delta^{\prime}=\left(1,-1\right)\allowbreak\left(3,-3\right)\left(5,-5\right)\left(7,-7\right)\left(9,-9\right)\), we calculate that \(\hat{\pi}=\left(1,-9\right)\left(-1,9\right)\left(2,7\right)\allowbreak\left(-2,-7\right)\left(3,4\right)\left(-3,-4\right)\left(5,10\right)\left(-5,-10\right)\left(6,-8\right)\left(-6,8\right)\), as expected.  If we have \(\gamma_{+}=\left(1,2,3,4,5,6\right)\left(7,8,9,10\right)\) and \(\gamma_{-}=\left(-1,-2,-3,-4,-5,-6\right)\allowbreak\left(-7,-8,-9,-10\right)\), we calculate that \(\gamma_{-}^{-1}\hat{\pi}\gamma_{+}=\left(1,7,-5,-9\right)\allowbreak\left(-1,9,5,-7\right)\allowbreak\left(2,4,10\right)\left(-2,-10,-4\right)\left(3\right)\left(-3\right)\left(6,-8\right)\left(-6,8\right)\).
\end{example}

\section{Proof of the Main Formula}
\label{proof}

We now present the proof of the main formula.  We will need a lemma to show the cycles of \(\gamma_{-}^{-1}\hat{\pi}\gamma_{+}\) do indeed come in front-back pairs as expected, so one of each pair may be discarded.  In the course of this proof we will show some properties of \(\hat{\pi}\) which will be useful later.  In the more general matrix model with arbitrary transposes (Section~\ref{arbitrary transposes}), we will use a permutation \(\delta^{\prime}\) which might reverse the signs on any \(k\) and \(-k\), so the following proofs will not depend on these \(k\) being odd.

\begin{lemma}
\label{redundant cycles}
With permutations as given in Notations~\ref{gamma}, \ref{delta}, \ref{pi hat} and \ref{more gammas}, if the cycle \(\left(k_{1},\ldots,k_{m}\right)\) appears in the cycle decomposition of \(\gamma_{-1}^{-1}\hat{\pi}\gamma_{+}\), then the cycle \(\left(-k_{m},\ldots,-k_{1}\right)\) also appears as a distinct cycle.

Furthermore, the permutation \(\hat{\pi}\) is a pairing with \(\left|\hat{\pi}\left(k\right)\right|=\pi\left(\left|k\right|\right)\).
\begin{proof}
We first show that the cycle \(\left(-k_{m},\ldots,-k_{1}\right)\) appears.  We do so by showing that \(\gamma_{-}^{-1}\hat{\pi}\gamma_{+}=\delta\left(\gamma_{-}^{-1}\hat{\pi}\gamma_{+}\right)^{-1}\delta\); that is, if we reverse all the cycles of \(\gamma_{-}^{-1}\hat{\pi}\gamma_{+}\) and change the signs on each integer, we get the same permutation back.

We can see that \(\delta\left(\gamma_{-}^{-1}\hat{\pi}\gamma_{+}\right)^{-1}\delta=\delta\gamma_{+}^{-1}\delta\delta\hat{\pi}^{-1}\delta\delta\gamma_{-}\delta\).  Since conjugation by \(\delta\) changes the signs on each integer in the cycle notation of a permutation, we can see that \(\delta\gamma_{+}^{-1}\delta=\gamma_{-}^{-1}\) and \(\delta\gamma_{-}\delta=\gamma_{+}\).

We can see that \(\hat{\pi}\) is self-inverse: \(\hat{\pi}=\left(\delta^{\prime}\pi\delta\pi\delta^{\prime}\right)^{-1}=\delta^{\prime}\pi\delta\pi\delta^{\prime}\), since \(\delta^{\prime}\), \(\pi\) and \(\delta\) are all involutions.

We now show that \(\hat{\pi}\) is unchanged by conjugation by \(\delta\).  \(\pi\delta\pi\) is \(\delta\) conjugated by \(\pi\), so its cycle decomposition is that of \(\delta\) with each positive integer \(k\) replaced with \(l:=\pi\left(k\right)\).  Its cycles are then of the form \(\left(l,-k\right)\), \(1\leq k\leq 2n\).

Conjugating \(\pi\delta\pi\) by \(\delta^{\prime}\) reverses the sign on some integers, so the sign on either or both integers in the cycle \(\left(l,-k\right)\) may be reversed.  In that case, the signs on their negatives in the cycle \(\left(-l,k\right)\) (which must also appear) will also be reversed.  So \(\hat{\pi}=\delta^{\prime}\pi\delta\pi\delta^{\prime}\) consists of cycles of the form \(\left(l,\pm k\right)\) and \(\left(-l,\mp k\right)\), \(1,\leq l\leq 2n\) (each cycle will be represented twice here).  Since \(k\neq l\) for all \(k\), each of these cycles contains two distinct elements, so \(\hat{\pi}\) is a pairing.  We have proven the stated properties of \(\hat{\pi}\).

Conjugating \(\hat{\pi}\) by \(\delta\) reverses the signs on all integers in the cycle decomposition of \(\hat{\pi}\), which changes cycles of the form \(\left(l,\pm k\right)\) into \(\left(-l,\mp k\right)\) and {\em vice versa}.  So we see that \(\hat{\pi}\) is unchanged by conjugation by \(\delta\).

We thus see that \(\delta\left(\gamma_{-}^{-1}\hat{\pi}\gamma_{+}\right)^{-1}\delta=\gamma_{-}^{-1}\hat{\pi}\gamma_{+}\).  We now show that these cycles must be distinct.  If not, for each \(i\) \(k_{i}\) and \(-k_{i}\) must occur in the same cycle.  Thus there must be an \(i\) and \(j\) with \(i<j\) such that \(k_{i}=-k_{j}\).  Then \(k_{i+t}=-k_{j-t}\) (by the above arguments) for all integers \(t\).  We can thus choose our \(i\) such that either \(k_{i}=-k_{i}\) or \(k_{i}=-k_{i+1}\).  Clearly the first case cannot occur.  In the second case, \(\gamma_{-}^{-1}\hat{\pi}\gamma_{+}\left(k_{i}\right)=k_{i+1}=-k_{i}\), so \(\hat{\pi}\gamma_{+}\left(k_{i}\right)=\gamma_{-}\left(-k_{i}\right)=\delta\gamma_{+}\delta\left(\delta\left(k_{i}\right)\right)=-\gamma_{+}\left(k_{i}\right)\).  However, as we saw above, the cycles of \(\hat{\pi}\) are of the form \(\left(\pm k,\pm\pi\left(k\right)\right)\) (with all four choices of \(+\) and \(-\) possible, depending on \(k\)).  Since \(\pi\) is a pairing, it never maps an integer to itself, so \(\hat{\pi}\) cannot map an an integer to another with the same absolute value.  So it is not possible for \(\hat{\pi}\left(\gamma_{+}\left(k_{i}\right)\right)\) to be equal to \(-\gamma_{+}\left(k_{i}\right)\).  Thus \(k_{i}\) and \(-k_{i}\) must appear in different cycles, so the cycles \(\left(k_{1},\ldots,k_{m}\right)\) and \(\left(-k_{m},\ldots,-k_{1}\right)\) must be distinct.
\end{proof}
\end{lemma}

We can discard the redundant cycles by any method, such as by taking only cycles whose integer of lowest absolute value is positive (called {\em particular} in \cite{MR2132270}).

\begin{notation}
\label{permutation}
We denote by \(\gamma_{-}\hat{\pi}\gamma_{+}/2\) the set of particular cycles.  We will also use this symbol to denote the set of integers appearing in these cycles, and we will think of the permutation as acting on this set.
\end{notation}

We now prove the main theorem:

\begin{theorem}
\begin{multline}
a_{r}\left(Y_{1},\ldots,Y_{r}\right)=\mathbb{E}\left(Y_{1}\cdots Y_{r}\right)\\=N^{-n-r}\sum_{\pi\in{\cal P}_{2}\left(2n\right)}\mathrm{Tr}_{\gamma_{-}^{-1}\hat{\pi}\gamma_{+}/2}\left(D_{1},\ldots,D_{2n}\right)
\label{main formula}
\end{multline}
where the notation is as defined above in Notations~\ref{gamma}, \ref{delta}, \ref{pi hat}, \ref{more gammas} and \ref{permutation}.

\begin{proof}
As shown in Lemma~\ref{lemma}:
\begin{multline*}
a_{r}\left(Y_{1},\ldots,Y_{r}\right)=N^{-n-r}\sum_{\pi\in{\cal P}_{2}\left(2n\right)}\sum_{\substack{i:\left[2n\right]\rightarrow\left[N\right]:i=i\circ\pi\nonumber\\j:\left[2n\right]\rightarrow\left[M\right]:j=j\circ\pi}}d^{\left(1\right)}_{i_{1}i_{\gamma\left(1\right)}}d^{\left(2\right)}_{j_{2}j_{\gamma\left(2\right)}}\cdots\\d^{\left(2n-1\right)}_{i_{2n-1}i_{\gamma\left(2n-1\right)}}d^{\left(2n\right)}_{j_{2n}j_{\gamma\left(2n\right)}}\textrm{.}
\end{multline*}

We will continue to denote \(D_{-k}:=D_{k}^{T}\).  Likewise, define \(d^{\left(-k\right)}_{ij}\) to be the \(ij\)th entry of \(D_{-k}\) (that is, the \(ji\)th entry of \(D_{k}\)).  Then \(d^{\left(-k\right)}_{ij}=d^{\left(k\right)}_{ji}\).

We note that the indices appearing on matrix \(D_{k}\) are \(\iota_{k}\) and \(\iota_{\gamma\left(k\right)}\), where \(\iota\) represents \(i\) or \(j\), depending on which indices appear on that particular matrix.  So the term from \(D_{k}\) which appears in our expression is \(d^{\left(k\right)}_{\iota_{k}\iota_{\gamma\left(k\right)}}\).  If we ignore the sign on the subscripts of indices (that is, we define \(\iota_{-k}:=\iota_{k}\)), this is equal to \(d^{\left(-k\right)}_{\iota_{\gamma_{-}\left(-k\right)}\iota_{-k}}\).  Since \(\gamma_{+}\) has no effect on negative integers and \(\gamma_{-}\) has no effect on positive integers, both of these expressions are equal to \(d^{\left(k\right)}_{\iota_{\gamma_{-}\left(k\right)}\iota_{\gamma_{+}\left(k\right)}}=d^{\left(-k\right)}_{\iota_{\gamma_{-}\left(-k\right)}\iota_{\gamma\left(-k\right)}}\).  We can replace each term \(d^{\left(k\right)}_{\iota_{k}\iota_{\gamma\left(k\right)}}\) by \(d^{\left(\pm k\right)}_{\iota_{\gamma_{-}\left(\pm k\right)}\iota_{\gamma_{+}\left(\pm k\right)}}\), where \(\pm k\) is whichever of \(k\) and \(-k\) appears in \(\gamma_{-}^{-1}\hat{\pi}\gamma_{+}/2\), and since exactly one of \(k\) and \(-k\) appears for each \(1\leq k\leq 2n\), we can choose the superscripts to be the elements of \(\gamma_{-}^{-1}\hat{\pi}\gamma_{+}/2\).

We have \(i=i\circ\pi\) and \(j=j\circ\pi\), and since \(\left|\hat{\pi}\left(k\right)\right|=\pi\left(\left|k\right|\right)\) (Lemma~\ref{redundant cycles}), \(\iota_{k}=\iota_{\hat{\pi}\left(k\right)}\) for all \(\iota\) and \(k\), we can replace the conditions \(i=i\circ\pi\) and \(j=j\circ\pi\) with \(i=i\circ\hat{\pi}\) and \(j=j\circ\hat{\pi}\).

We can now write our expression:

\[a_{r}\left(Y_{1},\ldots,Y_{r}\right)=N^{-n-r}\sum_{\pi\in{\cal P}_{2}\left(2n\right)}\sum_{\substack{i:\left[2n\right]\rightarrow\left[N\right]:i=i\circ\hat{\pi}\\j:\left[2n\right]\rightarrow\left[M\right]:j=j\circ\hat{\pi}}}\prod_{k\in\gamma_{-}^{-1}\hat{\pi}\gamma_{+}/2}d^{\left(k\right)}_{\iota_{\gamma_{-}\left(k\right)}\iota_{\gamma_{+}\left(k\right)}}\textrm{.}\]

We now show that each of these constraints \(\iota_{k}=\iota_{\hat{\pi}\left(k\right)}\) equates a first index with a second index in this expression.  A second index is of the form \(\iota_{\gamma_{+}\left(k\right)}\) for some \(k\in\gamma_{-}^{-1}\hat{\pi}\gamma_{+}/2\).  The index \(\iota_{\gamma_{+}\left(k\right)}\) is equated with \(\iota_{\hat{\pi}\gamma_{+}\left(k\right)}\).  We can see that \(\hat{\pi}\gamma_{+}\left(k\right)=\gamma_{-}\left(\gamma_{-}^{-1}\hat{\pi}\gamma_{+}\left(k\right)\right)\).  The integer \(\gamma_{-}^{-1}\hat{\pi}\gamma_{+}\left(k\right)\) is the next element after \(k\) in the cycle of \(\gamma_{-}^{-1}\hat{\pi}\gamma_{+}\), so it appears in the same cycle and hence must also be an element of \(\gamma_{-}^{-1}\hat{\pi}\gamma_{+}/2\).  Then \(\gamma_{-}\left(\gamma_{-}^{-1}\hat{\pi}\gamma_{+}\left(k\right)\right)\) is the subscript on the first index on the entry from \(D_{\gamma_{-}^{-1}\hat{\pi}\gamma_{+}\left(k\right)}\).  Thus the second indices are constrained by the values of the first indices, which are otherwise unconstrained (\(\hat{\pi}\) being a pairing).

Thus, instead of the constraints \(i=i\circ\hat{\pi}\) and \(j=j\circ\hat{\pi}\), we can rewrite all the second indices as \(\hat{\pi}\gamma_{+}\left(k\right)\).  This leaves us with terms of the form \(d^{\left(k\right)}_{\iota_{\gamma_{-}\left(k\right)}\iota_{\hat{\pi}\gamma_{+}\left(k\right)}}\).  We renumber our indices by letting \(\iota^{\prime}_{k}=\iota_{\gamma_{-}\left(k\right)}\) (legitimate since \(\gamma_{-}\) is a permutation).  Our expression is now:
\[a_{r}\left(Y_{1},\ldots,Y_{r}\right)=N^{-n-r}\sum_{\pi\in{\cal P}_{2}\left(2n\right)}\sum_{\substack{i:\left[2n\right]\rightarrow\left[N\right]\\j:\left[2n\right]\rightarrow\left[M\right]}}\prod_{k\in\gamma_{-}^{-1}\hat{\pi}\gamma_{+}/2}d^{\left(k\right)}_{\iota^{\prime}_{k}\iota^{\prime}_{\gamma_{-}\hat{\pi}\gamma_{+}\left(k\right)}}\textrm{.}\]
We recognize the term for each \(\hat{\pi}\) as the expression for the trace along the permutation \(\gamma_{-}^{-1}\hat{\pi}\gamma_{+}/2\) as given in Lemma~\ref{trace along permutation}.
\end{proof}
\end{theorem}

We now demonstrate that the cumulants \(k_{r}\left(Y_{1},\ldots,Y_{r}\right)\) are given by sums over the connected figures, or, in terms of the permutations, those where the subgroup of \(S_{2n}\) generated by \(\gamma\) and \(\pi\) acts transitively on \(\left[2n\right]\).  We demonstrate this by showing that the quantities defined in this manner satisfy the moment-cumulant formula (\ref{moment-cumulant formula}).

We note that the elements of a cycle of \(\gamma\) must be contained in the same orbit of \(\langle\gamma,\pi\rangle\).

\begin{definition}
We say that pairing \(\pi\) {\em connects} cycles \((\left(2\left(n_{1}+\cdots+n_{k-1}\right)\right.\allowbreak\left.+1,\ldots,2\left(n_{1}+\cdots+n_{k}\right)\right)\) and \(\left(2\left(n_{1}+\cdots+n_{l-1}\right)+1,\ldots,2\left(n_{1}+\cdots+n_{l}\right)\right)\) (which we will refer to as the \(i\)th and \(j\)th cycle, respectively) if their elements appear in the same orbit of \(\langle\gamma,\pi\rangle\).  If \(\langle\gamma,\pi\rangle\) acts transitively, we say that \(\pi\) connects the cycles of \(\gamma\).

We can construct a partition of \(\left[r\right]\) where \(i\) and \(j\) appear in the same block when the \(i\)th and \(j\)th cycles of \(\gamma\) are contained in the same orbit of \(\langle\gamma,\pi\rangle\).  We say that \(\pi\) {\em induces} this partition on \(\left[r\right]\).
\end{definition}

\begin{proposition}
\label{cumulants}
The cumulants of the \(Y_{k}\) are given by
\[k_{r}\left(Y_{1},\ldots,Y_{r}\right)=\sum_{\substack{\pi\in{\cal P}_{2}\left(2n\right)\\\textrm{\(\langle\gamma,\pi\rangle\) acts transitively}}}\mathrm{Tr}_{\gamma_{-}^{-1}\hat{\pi}\gamma_{+}/2}\left(D_{1},\ldots,D_{2n}\right)\]
\begin{proof}
Define \(r\)-linear functions \(\tilde{k}_{r}\) by
\[\tilde{k}_{r}\left(Y_{1},\ldots,Y_{r}\right):=\sum_{\substack{\pi\in{\cal P}_{2}\left(2n\right)\\\textrm{\(\langle\gamma,\pi\rangle\) acts transitively}}}\mathrm{Tr}_{\gamma_{-}^{-1}\hat{\pi}\gamma_{+}/2}\left(D_{1},\ldots,D_{2n}\right)\]
and for a partition \(\rho\in{\cal P}\left(r\right)\), let
\[\tilde{k}_{\rho}\left(Y_{1},\ldots,Y_{r}\right):=\prod_{V=\left\{i_{1},\ldots,i_{s}\right\}\in\rho}\tilde{k}_{s}\left(Y_{i_{1}},\ldots,Y_{i_{s}}\right)\textrm{.}\]

For each \(\rho\in{\cal P}\left(r\right)\), each term \(\tilde{k}_{s}\left(Y_{i_{1}},\ldots,Y_{i_{s}}\right)\) corresponding to block \(V=\left\{i_{1},\ldots,i_{s}\right\}\in\rho\) is a sum over all pairings on \(\left[2\left(n_{i_{1}}+\cdots+n_{i_{s}}\right)\right]\) connecting the cycles of permutation \(\left(1,\ldots,2n_{i_{1}}\right)\cdots\left(2\left(n_{i_{1}}+\cdots+n_{i_{s-1}}\right),\ldots,\right.\allowbreak\left.2\left(n_{i_{1}}+\cdots+n_{i_{s}}\right)\right)\).  If we identify the elements of \(\left[2\left(n_{i_{1}}+\cdots+n_{i_{s}}\right)\right]\) with the elements of the \(i_{1}\textrm{th},\ldots,i_{s}\textrm{th}\) cycles of \(\gamma\) (in order), we can see that this permutation is \(\gamma\) restricted to the elements of these cycles.

The term \(\tilde{k}_{\rho}\left(Y_{1},\ldots,Y_{r}\right)\) can thus be interpreted as a product of sums which, when expanded out, is the sum of terms, each corresponding to a choice for each block \(V\in\rho\) of a pairing connecting the cycles indexed by \(V\), that is, a pairing of \(\left[2n\right]\) inducing the partition \(\rho\).  Any pairing of \(\left[2n\right]\) inducing the partition \(\rho\) consists of such a choice.  Thus \(\tilde{k}_{\rho}\left(Y_{1},\ldots,Y_{s}\right)\) is the sum of all terms \(\mathrm{Tr}_{\gamma_{-}^{-1}\hat{\pi}\gamma_{+}/2}\left(D_{1},\ldots,D_{2n}\right)\) for \(\pi\) inducing \(\rho\).

Since any pairing \(\pi\) must induce some partition \(\rho\in{\cal P}\left(r\right)\), the moment \(a_{r}\left(Y_{1},\ldots,Y_{s}\right)\) is the sum of the \(\tilde{k}_{\rho}\left(Y_{1},\ldots,Y_{r}\right)\) for all \(\rho\in{\cal P}\left(r\right)\).  Thus the \(\tilde{k}_{s}\) satisfy the moment cumulant formula (\ref{moment-cumulant formula}).
\end{proof}
\end{proposition}

\section{Asymptotics}
\label{asymptotics}

For \(\sigma\in S_{n}\), let \(\#\left(\pi\right)\) denote the number of cycles of \(\sigma\) (that is, the number of orbits in its action on \(\left[2n\right]\), and for \(\sigma,\tau\in S_{n}\), let \(\#\left(\langle\sigma,\tau\rangle\right)\) denote the number of orbits of the subgroup \(\langle\sigma,\tau\rangle\).  In place of the Euler characteristic, we will use the inequality
\[\#\left(\sigma\right)+\#\left(\sigma^{-1}\tau\right)+\#\left(\tau\right)\leq n+2\#\left(\langle\sigma,\tau\rangle\right)\]
(see \cite{MR1396978, MR2052516}).

Since \(\gamma_{-}^{-1}\hat{\pi}\gamma_{+}\) is conjugate to \(\gamma_{-}\gamma_{-}^{-1}\hat{\pi}\gamma_{+}\gamma_{-}^{-1}=\hat{\pi}\gamma_{+}\gamma_{-}^{-1}\), these two permutations have the same number of cycles, so we find that
\[\#\left(\hat{\pi}\right)+\#\left(\gamma_{-}^{-1}\hat{\pi}\gamma_{+}\right)+\#\left(\gamma_{+}\gamma_{-}^{-1}\right)\leq 4n+2\#\left(\langle\gamma_{+}\gamma_{-}^{-1},\hat{\pi}\rangle\right).\]

We showed in Lemma~\ref{redundant cycles} that \(\hat{\pi}\) is a pairing, so \(\#\left(\hat{\pi}\right)=2n\).  We can calculate that \(\gamma_{+}\gamma_{-}^{-1}=\left(1\,\ldots,2n_{1}\right)\cdots\left(2\left(n_{1}+\cdots+n_{r-1}\right)+1,\ldots,2n\right)\left(-2n_{1},\right.\allowbreak\left.\ldots,-1\right)\cdots\left(-2n,\ldots,-2\left(n_{1}+\ldots+n_{r-1}\right)\right)\), which has \(2r\) cycles.  The subgroup \(\langle\gamma_{+}\gamma_{-}^{-1},\hat{\pi}\rangle\) cannot have more orbits than the subgroup generated by \(\gamma_{+}\gamma_{-}^{-1}\) alone, so it must have at most \(2r\) orbits.  The permutation \(\gamma_{-}^{-1}\hat{\pi}\gamma_{+}/2\) has half as many cycles as \(\gamma_{-}^{-1}\hat{\pi}\gamma_{+}\).  So we find
\[\#\left(\gamma_{-}^{-1}\hat{\pi}\gamma_{+}/2\right)\leq n+r.\]
If we express the moment \(a_{r}\left(Y_{1},\ldots,Y_{r}\right)\) in terms of the normalized trace, we find that the leading order terms (those for which the equality holds) have order \(N^{0}\).

When calculating the cumulants, we only consider pairings \(\pi\) such that \(\langle\gamma,\pi\rangle\) acts transitively.  Topologically, we can see that for such a \(\pi\), \(\langle\gamma_{+}\gamma_{-}^{-1},\hat{\pi}\rangle\) has at most \(2\) cycles.  Algebraically, for such a \(\pi\), we can find a sequence of cycles \(\left(k,l\right)\) of \(\pi\) connecting the \(s\)th and \(t\)th cycles of \(\gamma\) for any \(s,t\in\left[r\right]\).  As shown in Lemma~\ref{redundant cycles}, if cycle \(\left(k,l\right)\) appears in \(\pi\), cycles \(\left(k,\pm l\right)\) and \(\left(-k,\mp l\right)\) appear in \(\hat{\pi}\), so we can find a sequence of cycles of \(\hat{\pi}\) connecting the \(s\)th cycle of \(\gamma_{+}\) to the \(t\)th cycle of either \(\gamma_{+}\) or \(\gamma_{-}\) (or possibly both), and likewise with the \(s\)th cycle of \(\gamma_{-}\).  Thus the subgroup \(\langle\gamma_{+}\gamma_{-}^{-1},\hat{\pi}\rangle\) can then have at most \(2\) orbits.  Thus
\[\#\left(\gamma_{-}^{-1}\hat{\pi}\gamma_{+}/2\right)\leq n-r+2.\]
The order of \(k_{r}\left(Y_{1},\ldots,Y_{s}\right)\) is then \(N^{2-2r}\) and
\begin{multline}
\label{asymptotic value}
\lim_{N\rightarrow\infty}N^{2r-2}k_{r}\left(Y_{1},\ldots,Y_{r}\right)\\=\sum_{\substack{\pi\in{\cal P}_{2}\left(2n\right)\\\textrm{\(\langle\gamma,\pi\rangle\) acts transitively}}}\lim_{N\rightarrow\infty}\mathrm{tr}_{\gamma_{-}^{-1}\hat{\pi}\gamma_{+}/2}\left(D_{1},\ldots,D_{2n}\right)\textrm{.}
\end{multline}

\begin{theorem}
As \(N\rightarrow\infty\), the random variable \(Y_{k}\) converges weakly to the constant random variable with value \(\lim_{N\rightarrow\infty}k_{1}\left(Y_{k}\right)\).

Furthermore, the centred and rescaled random variables \(N\left(Y_{k}-\mathbb{E}\left(Y_{k}\right)\right)\) converge weakly to a centred multivariate Gaussian random variable with covariances \(\lim_{N\rightarrow\infty}N^{2}k_{2}\left(Y_{k},Y_{l}\right)\).

\begin{proof}
We will use the method of moments in both cases.  In each, convergence in moments implies weak convergence (see, e.g. \cite{MR1767078} for a proof of the method of moments and, e.g. \cite{MR1324786}, page 389, for a proof that a Gaussian distribution is determined by its moments).

We can see that the first cumulant of \(Y_{k}\) converges to the value given by (\ref{asymptotic value}) as \(N\rightarrow\infty\), while all higher cumulants converge to zero since they are of higher order in \(\frac{1}{N}\).  These are the cumulants of a constant random variable.

Centring the random variable \(Y_{k}\) does not change the value of any cumulant \(k_{r}\) for \(r>1\) (since the cumulants are multilinear and any set of constant random variables are independent and hence have zero higher cumulants).  After mutliplying by \(N\), the first cumulant is zero, the second is \(\lim_{N\rightarrow\infty}N^{2}k_{2}\left(Y_{k},Y_{l}\right)\) as given in (\ref{asymptotic value}), and all higher cumulants converge to zero as \(N\rightarrow\infty\).  These are the cumulants of a multivariate Gaussian.
\end{proof}
\end{theorem}

\begin{remark}
As in \cite{MR2052516}, highest order terms may be described in terms of planar diagrams, although we will not do so rigorously.  If each face is drawn in the plane in one of its orientations, then a pairing corresponding to a highest order term is one which may be drawn without crossings while the twistedness of each edge identification respects the relative orientations of the two faces.

The asymptotic expected value of a \(Y_{k}\) is the same as in the complex case, since any twisted edge identification will result in a non-orientable surface, which will not be a sphere.  The \(k\)th moment will have \(2^{r-1}\) times as many terms as in the complex case, since there are \(2^{r-1}\) possible relative orientations of the faces, but the bijection is not natural, so unless the constant matrices are chosen in a particularly symmetric way, the moment will not be \(2^{r-1}\) times the value in the complex case.
\end{remark}

\section{More General Matrix Models}
\label{generalisations}

\subsection{Several Wishart Matrices}

Instead of a single random matrix \(X\), we can consider a family of random matrices indexed by a Hilbert space \({\cal H}\): \(\left\{X_{g}:g\in{\cal H}\right\}\).  We let the \(ij\)th entry of \(X_{g}\) be \(\frac{1}{\sqrt{N}}f_{ij}^{\left(g\right)}\), where the \(f_{ij}^{\left(g\right)}\) are \(N\left(0,1\right)\) random variables which are independent for different choices of \(i\) and \(j\) and such that \(\mathbb{E}\left(f_{ij}^{\left(g\right)}f_{ij}^{\left(h\right)}\right)=\langle g,h\rangle\) for \(g,h\in{\cal H}\).

Applying the Wick formula creates the same constraints on nonzero terms for a given pairing \(\pi\).  However, the value of the expected value expression (\ref{expected value expression}), when it is nonzero, will now be
\[\prod_{\left\{k,l\right\}\in\pi}\mathbb{E}\left(f_{i_{k}j_{k}}^{\left(g_{k}\right)}f_{i_{l}j_{l}}^{\left(g_{l}\right)}\right)=\prod_{\left\{k,l\right\}\in\pi}\langle g_{k},g_{l}\rangle\]
The term in (\ref{sum over pairings}) associated to the pairing \(\pi\) will be weighted by this value, so in place of (\ref{main formula}), we get
\[a_{r}\left(Y_{1},\ldots,Y_{r}\right)=N^{-n-r}\sum_{\pi\in{\cal P}_{2}\left(2n\right)}\prod_{\left\{k,l\right\}\in\pi}\langle g_{k},g_{l}\rangle\mathrm{Tr}_{\gamma_{-}^{-1}\hat{\pi}\gamma_{+}/2}\left(D_{1},\ldots,D_{2n}\right)\textrm{.}\]

In particular, we can consider several independent Wishart matrices: if the pairing connects terms from different matrices, the resulting term will be zero.  We can colour the edges according to the different matrices, and consider only pairings which match colours.

\subsection{\(q\)-Wishart Matrices}

We can also consider matrices whose entries are \(q\)-commutative Gaussian random variables, that is, noncommutative random variables satisfying the \(q\)-Wick formula, as discussed in \cite{MR2439565}.

\begin{definition}
Let \(\pi\in{\cal P}_{2}\left(n\right)\) be a pairing.  We call a {\em crossing} a pair of blocks \(\left\{i,j\right\},\left\{k,l\right\}\in\pi\) such that \(i<k<j<l\).  We denote by \(\mathrm{cr}\left(\pi\right)\) the number of crossings in \(\pi\).
\end{definition}

\begin{definition}
Noncommutative random variables \(\left\{X_{g}\right\}_{g\in{\cal H}}\), for some Hilbert space \({\cal H}\) are said to be {\em real \(q\)-Gaussian} if they satisfy the \(q\)-Wick formula:
\[\mathbb{E}\left(X_{g_{1}}\cdots X_{g_{n}}\right)=\sum_{\pi\in{\cal P}_{2}\left(n\right)}q^{\mathrm{cr}\left(\pi\right)}\prod_{\left\{k,l\right\}\in\pi}\langle g_{k},g_{l}\rangle\]
\end{definition}

Again, the expected value expression (\ref{expected value expression}) can be replaced by this one, and the expression corresponding to pairing \(\pi\) becomes a weight on the term in (\ref{sum over pairings}) corresponding to \(\pi\), so in place of (\ref{main formula}), we get
\begin{multline*}
a_{r}\left(Y_{1},\ldots,Y_{r}\right)\\=N^{-n-r}\sum_{\pi\in{\cal P}_{2}\left(2n\right)}q^{\mathrm{cr}\left(\pi\right)}\prod_{\left\{k,l\right\}\in\pi}\langle g_{k},g_{l}\rangle\mathrm{Tr}_{\gamma_{-}^{-1}\hat{\pi}\gamma_{+}/2}\left(D_{1},\ldots,D_{2n}\right)\textrm{.}
\end{multline*}

\begin{remark}
In an expression with only one trace, crossings correspond to the geometric notion of crossings.  However, in an expression with several traces, the number of crossings depends on the order of the terms inside each trace, which is not expressed geometrically.
\end{remark}

\subsection{Arbitrary Transposes}
\label{arbitrary transposes}

We can also generalize our formula to expressions in which \(X\) and \(X^{T}\) terms may appear abitrarily, rather than alternatingly.  We let \(\varepsilon:\left[n\right]\rightarrow\left\{-1,1\right\}\) be a function encoding the positions of the transposes.  We denote by \(X^{\left(1\right)}\) the matrix \(X\) itself and by \(X^{\left(-1\right)}\) the transpose of \(X\).  In place of formula~(\ref{trace random variable}), we can calculate the moments and cumulants of random variables of the form
\[Y_{k}=\mathrm{tr}\left(X^{\left(\varepsilon\left(n_{1}+\cdots n_{k-1}+1\right)\right)}D_{n_{1}+\cdots+n_{k-1}}\cdots X^{\left(\varepsilon\left(n_{1}+\cdots n_{k}\right)\right)}D_{n_{1}+\cdots+n_{k}}\right)\]
where the \(D_{k}\) matrices are of the appropriate size for the above matrix multiplication to be defined.  We again let \(n=n_{1}+\cdots+n_{r}\).  It is no longer necessary that any of the \(n_{k}\) be even, although if \(n\) is odd there will be no pairings on \(\left[n\right]\), so the moments will all be equal to zero.

Instead of the indices \(i_{k}\) and \(j_{k}\), we sum over all indices \(1\leq\iota^{+}_{1},\ldots,\iota^{+}_{n}\leq M\) and \(1\leq\iota^{-}_{1},\ldots,\iota^{-}_{n}\leq N\) such that it is the \(\iota^{+}_{k}\iota^{-}_{k}\)th entry of \(X\) that is taken from its \(k\)th occurrence.  We then find that it is the \(\iota^{-\varepsilon\left(k\right)}_{k}\iota^{\varepsilon\left(\gamma\left(k\right)\right)}_{\gamma\left(k\right)}\)th entry of \(D_{k}\) which appears, that is,
\[d^{\left(\pm k\right)}_{\iota^{-\mathrm{sgn}\left(\pm k\right)\varepsilon\left(\gamma_{-}\left(\pm k\right)\right)}_{\gamma_{-}\left(\pm k\right)}\iota^{\mathrm{sgn}\left(\pm k\right)\varepsilon\left(\gamma_{+}\left(\pm k\right)\right)}_{\gamma^{+}\left(\pm k\right)}}.\]

We now define \(\delta^{\prime}\left(k\right):=\varepsilon\left(\left|k\right|\right)k\), so it reverses the sign on integer \(k\) when the \(\left|k\right|\)th occurrence of \(X\) appears with a transpose, and we construct \(\hat{\pi}\) using this \(\delta^{\prime}\).  Permutations \(\gamma\), \(\gamma_{-}\) and \(\gamma_{+}\) are defined as above.

The rest of the argument proceeds as above, and we find that
\[a_{r}\left(Y_{1},\ldots,Y_{r}\right)=N^{-\frac{n}{2}-r}\sum_{\pi\in{\cal P}_{2}\left(n\right)}\mathrm{Tr}_{\gamma_{-}^{-1}\hat{\pi}\gamma_{+}/2}\left(D_{1},\ldots,D_{n}\right)\textrm{.}\]

\begin{remark}
Since a real Wigner matrix can be written \(\frac{1}{2}\left(X+X^{T}\right)\), where \(X\) is a square matrix as defined in Definition~\ref{Wishart matrix}, we can expand to find that the \(n\)th moment of a compound Wigner matrix can be expressed:
\[\frac{1}{2^{n}}N^{-\frac{n}{2}-r}\sum_{\varepsilon:\left[n\right]\rightarrow\left\{-1,1\right\}}\sum_{\pi\in{\cal P}_{2}\left(n\right)}\mathrm{Tr}_{\gamma_{-}^{-1}\hat{\pi}\gamma_{+}/2}\left(D_{1},\ldots,D_{n}\right)\textrm{.}\]
Geometrically, this corresponds to a sum over all possible surface gluings with both directions of edge identification possible.
\end{remark}

\begin{example}
If we wish to calculate the quantity \(\mathbb{E}\left(\mathrm{Tr}\left(XD_{1}XD_{2}\right)\right)\), we use \(\delta^{\prime}=e\).  There is only one pairing on \(\left[2\right]\): \(\pi=\left(1,2\right)\), and \(\hat{\pi}=\left(1,-2\right)\left(2,-1\right)\).  We calculate that \(\gamma_{-}^{-1}\hat{\pi}\gamma_{+}=\left(1,-2\right)\left(2,-1\right)\).  Only the cycle \(\left(1,-2\right)\) is particular, so \(\gamma_{-}^{-1}\hat{\pi}\gamma_{+}/2=\left(1,-2\right)\), a permutation on \(\left\{1,-2\right\}\).  Thus
\[\mathbb{E}\left(\mathrm{Tr}\left(XD_{1}XD_{2}\right)\right)=N^{-2}\mathrm{Tr}\left(D_{1}D_{2}^{T}\right)=N^{-1}\mathrm{tr}\left(D_{1}D_{2}^{T}\right)\textrm{.}\]

If we let \(Z:=\frac{1}{2}\left(X+X^{T}\right)\) be a real Wigner matrix, then we find that
\begin{multline*}
\mathbb{E}\left(\mathrm{Tr}\left(ZD_{1}ZD_{2}\right)\right)=\frac{1}{4}\mathbb{E}\left(\mathrm{Tr}\left(XD_{1}XD_{2}\right)\right)+\frac{1}{4}\mathbb{E}\left(\mathrm{Tr}\left(XD_{1}X^{T}D_{2}\right)\right)+\\\frac{1}{4}\mathbb{E}\left(\mathrm{Tr}\left(X^{T}D_{1}XD_{2}\right)\right)+\frac{1}{4}\mathbb{E}\left(\mathrm{Tr}\left(X^{T}D_{1}X^{T}D_{2}\right)\right)\textrm{.}
\end{multline*}
Calculating the last three terms by the same methods, we find that
\begin{eqnarray*}
\mathbb{E}\left(\mathrm{Tr}\left(ZD_{1}ZD_{2}\right)\right)&=&\frac{1}{2}N^{-2}\mathrm{Tr}\left(D_{1}\right)\mathrm{Tr}\left(D_{2}\right)+\frac{1}{2}N^{-2}\mathrm{Tr}\left(D_{1}D_{2}^{T}\right)\\
&=&\frac{1}{2}\mathrm{tr}\left(D_{1}\right)\mathrm{tr}\left(D_{2}\right)+\frac{1}{2}N^{-1}\mathrm{tr}\left(D_{1}D_{2}^{T}\right)\textrm{.}
\end{eqnarray*}
Two sphere terms contribute to the first term on the right hand side, and two projective plane terms contribute to the second term.
\end{example}

\bibliography{paper}
\bibliographystyle{plain}

\end{document}